\documentclass[twoside,12pt]{article}
\usepackage{amssymb}
\usepackage{amsthm}
\usepackage{amsmath}
\usepackage{a4wide}
\usepackage[all]{xy}
\usepackage{fancyhdr}

\newtheorem{theorem}{Theorem}[section]
\newtheorem{proposition}[theorem]{Proposition}

\newtheorem{lemma}[theorem]{Lemma}
\theoremstyle{definition}
\newtheorem*{notation}{Notation}
\newtheorem*{Beweis}{Proof}
\newtheorem{definition}[theorem]{Definition}
\newtheorem{punto}[theorem]{}
\theoremstyle{remark}
\newtheorem{remark}[theorem]{Remark}
\newtheorem{ex}[theorem]{Example}
\newtheorem{c-ex}[theorem]{Counterexample}
\newtheorem{c-exs}[theorem]{Counterexamples}
\newtheorem{remarks}[theorem]{Remarks}
\CompileMatrices  \swapnumbers
\CompileMatrices
\input{tcilatex}

\begin{document}

\title{A Zariski Topology for Bicomodules and Corings\thanks{%
MSC (2000): 16W30, 16N60, 16D80\newline
Keywords: fully coprime (bi)comodules, fully cosemiprime (bi)comodules,
fully coprime corings, fully cosemiprime corings, fully coprime spectrum,
fully coprime coradical, Zariski topology, top (bi)comodules.}}
\author{\textbf{Jawad Y. Abuhlail}\thanks{%
Supported by King Fahd University of Petroleum $\&$ Minerals, Research
Project $\#\ $INT/296} \\
Department of Mathematical Sciences, Box $\#\ $5046\\
King Fahd University of Petroleum $\&$ Minerals\\
31261 Dhahran - Saudi Arabia\\
abuhlail@kfupm.edu.sa}
\date{}
\maketitle

\begin{abstract}
In this paper we introduce and investigate \emph{top }(\emph{bi})\emph{%
comodules} of corings, that can be considered as dual to top (bi)modules of
rings. The \emph{fully coprime spectra} of such (bi)comodules attains a 
\emph{Zariski topology}, defined in a way dual to that of defining the
Zariski topology on the prime spectra of (commutative) rings. We restrict
our attention in this paper to duo (bi)comodules (satisfying suitable
conditions) and study the interplay between the coalgebraic properties of
such (bi)comodules and the introduced Zariski topology. In particular, we
apply our results to introduce a Zariski topology on the fully coprime
spectrum of a given non-zero coring considered canonically as duo object in
its category of bicomodules.
\end{abstract}

\section{Introduction}

\qquad Several papers considered the so called \emph{top modules,} i.e.
modules (over commutative rings) whose spectrum of \emph{prime submodules}
attains a Zariski topology, e.g. \cite{Lu1999}, \cite{MMS1997}, \cite%
{Zha1999}. Dually, we introduce and investigate \emph{top }(\emph{bi})\emph{%
comodules} for corings and study their properties (restricting our attention
in this first paper to duo (bi)comodules satisfying suitable conditions). In
particular, we extend results of \cite{NT2001} on the topology defined on
the spectrum of (\emph{fully})\emph{\ coprime subcoalgebras} of a given
coalgebra over a base field to the general situation of a topology on the 
\emph{fully coprime spectrum} of a given non-zero bicomodule over a given
pair of non-zero corings.

Throughout, $R$ is a commutative ring with $1_{R}\neq 0_{R}$ and $A,B$ are $%
R $-algebras. With \emph{locally projective} modules, we mean those in the
sense of \cite{Z-H1976} (see also \cite{Abu2006}). We denote by $\mathcal{C}%
=(\mathcal{C},\Delta _{\mathcal{C}},\varepsilon _{\mathcal{C}})$ a non-zero $%
A$-coring with $_{A}\mathcal{C}$ flat and by $\mathcal{D}=(\mathcal{D}%
,\Delta _{\mathcal{D}},\varepsilon _{\mathcal{D}})$ a non-zero $B$-coring
with $\mathcal{D}_{B}$ flat, so that the categories $^{\mathcal{D}}\mathbb{M}%
^{\mathcal{C}}$ of $(\mathcal{D},\mathcal{C})$-bicomodules, $\mathbb{M}^{%
\mathcal{C}}$ of right $\mathcal{C}$-comodules and $^{\mathcal{D}}\mathbb{M}$
of left $\mathcal{D}$-comodules are Grothendieck.

After this brief introduction, we include in the second section some
preliminaries and extend some of our results in \cite{Abu2006} on fully
coprimeness in categories of comodules to fully coprimeness in categories of
bicomodules.

In the third and main section, we introduce a Zariski topology for
bicomodules. Let $M$ be a given non-zero $(\mathcal{D},\mathcal{C})$%
-bicomodule and consider the \emph{fully\ coprime spectrum}%
\begin{equation*}
\mathrm{CPSpec}(M):=\{K\mid K\subseteq M\text{ is a \emph{fully }}M\text{%
\emph{-coprime} }(\mathcal{D},\mathcal{C})\text{-subbicomodule}\}.
\end{equation*}%
For every $(\mathcal{D},\mathcal{C})$-subbicomodule $L\subseteq M,$ set $%
\mathcal{V}_{L}:=\{K\in \mathrm{CPSpec}(M)\mid K\subseteq L\}$ and $\mathcal{%
X}_{L}:=\{K\in \mathrm{CPSpec}(M)\mid K\nsubseteqq L\}.$ As in the case of
the spectra of prime submodules of modules over (commutative) rings (e.g. 
\cite{Lu1999}, \cite{MMS1997}, \cite{Zha1999}), the class of \emph{varieties}
$\xi (M):=\{\mathcal{V}_{L}\mid L\subseteq M$ is a $(\mathcal{D},\mathcal{C}%
) $-subbicomodule$\}$ satisfies all axioms of closed sets in a topological
space with the exception that $\xi (M)$ is \emph{not necessarily} closed
under finite unions. We say $M$ is a \emph{top bicomodule}, iff $\xi (M)$ is
closed under finite unions, equivalently iff $\tau _{M}:=\{\mathcal{X}%
_{L}\mid L\subseteq M$ is a $(\mathcal{D},\mathcal{C})$-subbicomodule$\}$ is
a topology (in this case we call $\mathbf{Z}_{M}:=(\mathrm{CPSpec}(M),\tau
_{M})$ a \emph{Zariski topology} of $M$). We then restrict our attention to
the case in which $M$ is a \emph{duo bicomodule} (i.e. every subbicomodule
of $M$ is fully invariant) satisfying suitable conditions. For such a
bicomodule $M$ we study the interplay between the coalgebraic properties of $%
M$ and the topological properties of $\mathbf{Z}_{M}.$

In the fourth section we give some applications and examples. Our main
application will be to non-zero corings which turn out to be \emph{duo }%
bicomodules in the canonical way. We also give some concrete examples that
establish some of the results in section three.

It is worth mentioning that several properties of the Zariski topology for
bicomodules and corings are, as one may expect, dual to those of the
classical Zariski topology on the prime spectrum of commutative rings (e.g. 
\cite{AM1969}, \cite{Bou1998}).

This paper is a continuation of \cite{Abu2006}. The ideas in both papers can
be transformed to investigate the notion of coprimeness in the sense of
Annin \cite{Ann2002} in categories of (bi)comodules and define a Zariski
topology on the spectrum of \emph{coprime sub}(\emph{bi})\emph{comodules} of
a given (bi)comodule. Moreover, different notions of primeness and
coprimeness in these papers can be investigated in categories of (bi)modules
over rings, which can be seen as bicomodules over the ground rings that
should be considered with the trivial coring structures (a different
approach has been taken in the recent work \cite{Wij2006}, where several
primeness and coprimeness conditions are studied in categories of modules
and then applied to categories of comodules of \emph{locally projective}
coalgebras over commutative rings). More generally, such (co)primeness
notions can be developed in more general Grothendieck categories. These and
other applications will be considered in forthcoming papers.

\section{Preliminaries}

\qquad All rings and their modules in this paper are assumed to be unital.
For a ring $T,$ we denote with $Z(T)$ the \emph{center} of $T$ and with $%
T^{op}$ the opposite ring of $T.$ For basic definitions and results on
corings and comodules, the reader is referred to \cite{BW2003}. A reference
for the topological terminology and other results we use could be any
standard book in general topology (notice that in our case, a \emph{compact
space} is not necessarily Hausdorff; such spaces are called \emph{%
quasi-compact} by some authors, e.g. \cite[I.9.1.]{Bou1966}).

\begin{punto}
\label{An=per}(e.g. \cite[17.8.]{BW2003}) For any $A$-coring $\mathcal{C},$
the dual module $^{\ast }\mathcal{C}:=\mathrm{Hom}_{A-}(\mathcal{C},A)$
(resp. $\mathcal{C}^{\ast }:=\mathrm{Hom}_{-A}(\mathcal{C},A)$) is an $%
A^{op} $-rings with unity $\varepsilon _{\mathcal{C}}$ and multiplication%
\begin{equation*}
(f\ast ^{l}g)(c):=\dsum f(c_{1}g(c_{2}))\text{ (resp. }(f\ast
^{r}g)(c):=\dsum g(f(c_{1})c_{2})\text{).}
\end{equation*}
\end{punto}

\begin{punto}
Let $M$ be a non-zero $(\mathcal{D},\mathcal{C})$-bicomodule. Then $M$ is a $%
(^{\ast }\mathcal{C},\mathcal{D}^{\ast })$-bimodule with actions%
\begin{equation*}
f\rightharpoonup m:=\dsum m_{<0>}f(m_{<1>})\text{ and }m\leftharpoonup
g:=\dsum g(m_{<-1>})m_{<0>},\text{ }f\in \text{ }^{\ast }\mathcal{C},\text{ }%
g\in \mathcal{D}^{\ast },\text{ }m\in M.
\end{equation*}%
Moreover, the set $^{\mathcal{D}}\mathrm{E}_{M}^{\mathcal{C}}:=$ $^{\mathcal{%
D}}\mathrm{End}^{\mathcal{C}}(M)^{op}$ of $(\mathcal{D},\mathcal{C})$%
-bicolinear endomorphisms of $M$ is a ring with multiplication the \emph{%
opposite} composition of maps, so that $M$ is canonically a $(^{\ast }%
\mathcal{C}\otimes _{R}\mathcal{D}^{\ast op},$ $^{\mathcal{D}}\mathrm{E}%
_{M}^{\mathcal{C}})$-bimodule. A $(\mathcal{D},\mathcal{C})$-subbicomodule $%
L\subseteq M$ is called \emph{fully invariant, }iff it is a right $^{%
\mathcal{D}}\mathrm{E}_{M}^{\mathcal{C}}$-submodule as well. We call $M\in $ 
$^{\mathcal{D}}\mathbb{M}^{\mathcal{C}}$ \emph{duo} (\emph{quasi-duo}), iff
every (simple) $(\mathcal{D},\mathcal{C})$-subbicomodule of $M$ is fully
invariant. If $_{A}\mathcal{C}$ and $\mathcal{D}_{B}$ are locally
projective, then $^{\mathcal{D}}\mathbb{M}^{\mathcal{C}}\simeq $ $^{\mathcal{%
D}}\mathrm{R}\mathrm{at}^{\mathcal{C}}(_{(\mathcal{D}^{\ast })^{op}}\mathbb{M%
}_{(^{\ast }\mathcal{C})^{op}})=$ $^{\mathcal{D}}\mathrm{R}\mathrm{at}^{%
\mathcal{C}}(_{^{\ast }\mathcal{C}}\mathbb{M}_{\mathcal{D}^{\ast }})$ (the
category of $(\mathcal{D},\mathcal{C})$\emph{-birational} $(^{\ast }\mathcal{%
C},\mathcal{D}^{\ast })$\emph{-bimodules}, e.g. \cite[Theorem 2.17.]{Abu2003}%
).
\end{punto}

\begin{notation}
Let $M$ be a $(\mathcal{D},\mathcal{C})$-bicomodule. With $\mathcal{L}(M)$
(resp. $\mathcal{L}_{f.i.}(M)$) we denote the lattice of (fully invariant) $(%
\mathcal{D},\mathcal{C})$-subbicomodules of $M$ and with $\mathcal{I}_{r}(^{%
\mathcal{D}}\mathrm{E}_{M}^{\mathcal{C}})$ (resp. $\mathcal{I}(^{\mathcal{D}}%
\mathrm{E}_{M}^{\mathcal{C}})$) the lattice of right (two-sided) ideals of $%
^{\mathcal{D}}\mathrm{E}_{M}^{\mathcal{C}}.$ With $\mathcal{I}_{r}^{f.g.}(^{%
\mathcal{D}}\mathrm{E}_{M}^{\mathcal{C}})\subseteq \mathcal{I}_{r}(^{%
\mathcal{D}}\mathrm{E}_{M}^{\mathcal{C}})$ (resp. $\mathcal{L}%
^{f.g.}(M)\subseteq \mathcal{L}(M)$) we denote the subclass of finitely
generated right ideals of $^{\mathcal{D}}\mathrm{E}_{M}^{\mathcal{C}}$ (the
subclass of $(\mathcal{D},\mathcal{C})$-subbicomodules of $M$ which are
finitely generated as $(B,A)$-bimodules). For $\varnothing \neq K\subseteq M$
and $\varnothing \neq I\subseteq $ $^{\mathcal{D}}\mathrm{E}_{M}^{\mathcal{C}%
}$ we set%
\begin{equation*}
\mathrm{An}(K):=\{f\in \text{ }^{\mathcal{D}}\mathrm{E}_{M}^{\mathcal{C}}|%
\text{ }f(K)=0\}\text{ and }\mathrm{Ke}(I):=\{m\in M\mid f(m)=0\text{ for
all }f\in I\}.
\end{equation*}
\end{notation}

\qquad In what follows we introduce some notions for an object in $^{%
\mathcal{D}}\mathbb{M}^{\mathcal{C}}:$

\begin{definition}
We say that a non-zero $(\mathcal{D},\mathcal{C})$-bicomodule $M$ is

\emph{self-injective, }iff for every $(\mathcal{D},\mathcal{C})$%
-subbicomodule $K\subseteq M,$ every $f\in $ $^{\mathcal{D}}\mathrm{Hom}^{%
\mathcal{C}}(K,M)$ extends to some $(\mathcal{D},\mathcal{C})$-bicolinear
endomorphism $\widetilde{f}\in $ $^{\mathcal{D}}\mathrm{E}_{M}^{\mathcal{C}%
}; $

\emph{self-cogenerator,} iff $M$ cogenerates $M/K$ in $^{\mathcal{D}}\mathbb{%
M}^{\mathcal{C}}$ $\forall $ $(\mathcal{D},\mathcal{C})$-subbicomodule $%
K\subseteq M;$

\emph{intrinsically injective}, iff $\mathrm{AnKe}(I)=I$ for every finitely
generated right ideal $I\vartriangleleft _{r}$ $^{\mathcal{D}}\mathrm{E}%
_{M}^{\mathcal{C}}.$

\emph{simple}, iff $M$ has no non-trivial $(\mathcal{D},\mathcal{C})$%
-subbicomodules;

\emph{subdirectly irreducible}\footnote{\emph{Subdirectly irreducible}
comodules were called \emph{irreducible} in \cite{Abu2006}. However, we
observed that such a terminology may cause confusion, so we choose to change
it in this paper to be consistent with the terminology used for modules
(e.g. \cite[9.11., 14.8.]{Wis1991}).}, iff $M$ contains a unique simple $(%
\mathcal{D},\mathcal{C})$-subbicomodule that is contained in every non-zero $%
(\mathcal{D},\mathcal{C})$-subbicomodule of $M$ (equivalently, iff $%
\dbigcap_{0\neq K\in \mathcal{L}(M)}\neq 0$).

\emph{semisimple}$\emph{,}$ iff $M=\mathrm{Corad}(M),$ where $\mathrm{Corad}%
(M):=\sum \{K\subseteq M\mid K$ is a simple $(\mathcal{D},\mathcal{C})$%
-subbicomodule$\}$ ($:=0,$ if $M$ has no simple $(\mathcal{D},\mathcal{C})$%
-subbicomodules).
\end{definition}

\begin{notation}
Let $M$ be a non-zero $(\mathcal{D},\mathcal{C})$-bicomodule. We denote with 
$\mathcal{S}(M)$ ($\mathcal{S}_{f.i.}(M)$) the class of simple $(\mathcal{D},%
\mathcal{C})$-subbicomodules of $M$ (non-zero fully invariant $(\mathcal{D},%
\mathcal{C})$-subbicomodules of $M$ with no non-trivial fully invariant $(%
\mathcal{D},\mathcal{C})$-subbicomodules). Moreover, we denote with $\mathrm{%
Max}_{r}(^{\mathcal{D}}\mathrm{E}_{M}^{\mathcal{C}})$ ($\mathrm{Max}(^{%
\mathcal{D}}\mathrm{E}_{M}^{\mathcal{C}})$) the class of maximal right
(two-sided) ideals of $^{\mathcal{D}}\mathrm{E}_{M}^{\mathcal{C}}.$ The
Jacobson radical (prime radical) of $^{\mathcal{D}}\mathrm{E}_{M}^{\mathcal{C%
}}$ is denoted by $\mathrm{Jac}(^{\mathcal{D}}\mathrm{E}_{M}^{\mathcal{C}})$
($\mathrm{Prad}(^{\mathcal{D}}\mathrm{E}_{M}^{\mathcal{C}})$).
\end{notation}

\begin{punto}
\label{ft=2}Let $M$ be a non-zero $(\mathcal{D},\mathcal{C})$-bicomodule. We 
$M$ has \emph{Property} $\mathbf{S}$ ($\mathbf{S}_{f.i.}$), iff $\mathcal{S}%
(L)\neq \varnothing $ ($\mathcal{S}_{f.i.}(L)\neq \varnothing $) for every
(fully invariant) non-zero $(\mathcal{D},\mathcal{C})$-subbicomodule $0\neq
L\subseteq M.$ Notice that if $M$ has $\mathbf{S},$ then $M$ is subdirectly
irreducible if and only if $L_{1}\cap L_{2}\neq 0$ for any two non-zero $(%
\mathcal{D},\mathcal{C})$-subbicomodules $0\neq L_{1},L_{2}\subseteq M.$
\end{punto}

\begin{lemma}
\label{simple}Let $M$ be a non-zero $(\mathcal{D},\mathcal{C})$-bicomodule.
If $B\otimes _{R}A^{op}$ is left perfect and $_{A}\mathcal{C},$ $\mathcal{D}%
_{B}$ are locally projective, then

\begin{enumerate}
\item every finite subset of $M$ is contained in a $(\mathcal{D},\mathcal{C}%
) $-subbicomodule $L\subseteq M$ that is finitely generated as a $(B,A)$%
-bimodule.

\item every non-zero $(\mathcal{D},\mathcal{C})$-subbicomodule $0\neq
L\subseteq M$ has a simple $(\mathcal{D},\mathcal{C})$-subbicomodule, so
that $M$ has Property $\mathbf{S.}$ If moreover, $M$ is quasi-duo, then $M$
has Property $\mathbf{S}_{f.i.}.$

\item $\mathrm{Corad}(M)\subseteq ^{e}M$ \emph{(}an essential $(\mathcal{D},%
\mathcal{C})$-subbicomodule\emph{)}.
\end{enumerate}
\end{lemma}

\begin{Beweis}
\begin{enumerate}
\item It's enough to show the assertion for a single element $m\in M.$ Let $%
\varrho _{M}^{\mathcal{C}}(m)=\dsum\limits_{i=1}^{n}m_{i}\otimes _{A}c_{i}$
and $\varrho _{M}^{\mathcal{D}}(m_{i})=\dsum\limits_{j=1}^{k_{i}}d_{i,j}%
\otimes _{B}m_{ij}$ for each $i=1,...,n.$ Since $_{A}\mathcal{C},$ $\mathcal{%
D}_{B}$ are locally projective, the $(^{\ast }\mathcal{C},\mathcal{D}^{\ast
})$-subbimodule $L:=$ $^{\ast }\mathcal{C}\rightharpoonup m\leftharpoonup 
\mathcal{D}^{\ast }\subseteq M$ is by \cite[Theorem 2.17.]{Abu2003} a $(%
\mathcal{D},\mathcal{C})$-subbicomodule. Moreover, $\{m_{i,j}\mid 1\leq
i\leq n,$ $1\leq j\leq k_{i}\}$ generates $_{B}L_{A},$ since%
\begin{equation*}
f\rightharpoonup m\leftharpoonup
g=[\dsum\limits_{i=1}^{n}m_{i}f(c_{i})]\leftharpoonup
g=\dsum\limits_{i=1}^{n}\dsum_{j=1}^{k_{i}}g(d_{i,j})m_{i,j}f(c_{i})\text{ }%
\forall \text{ }f\in \text{ }^{\ast }\mathcal{C}\text{ and }g\in \mathcal{D}%
^{\ast }.
\end{equation*}

\item Suppose $0\neq L\subseteq M$ is a $(\mathcal{D},\mathcal{C})$%
-subbicomodule with no simple $(\mathcal{D},\mathcal{C})$-subbicomodules. By
\textquotedblleft 1\textquotedblright , $L$ contains a non-zero $(\mathcal{D}%
,\mathcal{C})$-subbicomodule $0\neq L_{1}\subsetneqq M$ that is finitely
generated as a $(B,A)$-bimodule. Since $L$ contains no simple $(\mathcal{D},%
\mathcal{C})$-subbicomodules, for every $n\in \mathbb{N}$ we can pick (by
induction) a non-zero $(\mathcal{D},\mathcal{C})$-subbicomodule $0\neq
L_{n+1}\subsetneqq L_{n}$ that is finitely generated as a $B\otimes
_{R}A^{^{op}}$-module. In this way we obtain an infinite chain $%
L_{1}\supsetneqq L_{2}\supsetneqq ...\supsetneqq L_{n}\supsetneqq
L_{n+1}\supsetneqq ....$ of finitely generated $B\otimes _{R}A^{^{op}}$%
-submodules of $L$ (a contradiction to the assumption that $B\otimes
_{R}A^{^{op}}$ is left perfect, see \cite[Theorem 22.29]{Fai1976}).
Consequently, $L$ should contain at least one simple $(\mathcal{D},\mathcal{C%
})$-subbicomodule. Hence $M$ has property $\mathbf{S}.$ The last statement
is obvious.

\item For every non-zero $(\mathcal{D},\mathcal{C})$-subbicomodule $0\neq
L\subseteq M,$ we have by \textquotedblleft 1\textquotedblright\ $L\cap 
\mathrm{Corad}(M)=\mathrm{Corad}(L)\neq 0,$ hence $\mathrm{Corad}%
(M)\subseteq ^{e}M.\blacksquare $
\end{enumerate}
\end{Beweis}

Given a non-zero $(\mathcal{D},\mathcal{C})$-bicomodule $M,$ we have the
following annihilator conditions. The proofs are similar to the
corresponding results in \cite[28.1.]{Wis1991}, hence omitted:

\begin{punto}
\label{M-An-Ke}Let $M$ be a non-zero $(\mathcal{D},\mathcal{C})$-bicomodule
and consider the order-reversing mappings%
\begin{equation}
\mathrm{An}(-):\mathcal{L}(M)\rightarrow \mathcal{I}_{r}(^{\mathcal{D}}%
\mathrm{E}_{M}^{\mathcal{C}})\text{ and }\mathrm{Ke}(-):\mathcal{I}_{r}(^{%
\mathcal{D}}\mathrm{E}_{M}^{\mathcal{C}})\rightarrow \mathcal{L}(M).
\label{An-map-Ke}
\end{equation}

\begin{enumerate}
\item $\mathrm{An}(-)$ and $\mathrm{Ke}(-)$ restrict to order-reversing
mappings%
\begin{equation}
\mathrm{An}(-):\mathcal{L}_{f.i.}(M)\rightarrow \mathcal{I}(^{\mathcal{D}}%
\mathrm{E}_{M}^{\mathcal{C}})\text{ and }\mathrm{Ke}(-):\mathcal{I}(^{%
\mathcal{D}}\mathrm{E}_{M}^{\mathcal{C}})\rightarrow \mathcal{L}_{f.i.}(M).
\label{An-fi}
\end{equation}

\item For a $(\mathcal{D},\mathcal{C})$-subbicomodule $K\subseteq M:\mathrm{%
Ke}(\mathrm{An}(K))=K$ if and only if $M/K$ is $M$-cogenerated. So, if $M$
is self-cogenerator, then the map $\mathrm{An}(-)$ in (\ref{An-map-Ke}) and
its restriction in (\ref{An-fi}) are injective.

\item If $M$ is self-injective, then

\begin{enumerate}
\item $\mathrm{An}(\bigcap\limits_{i=1}^{n}K_{i})=\sum\limits_{i=1}^{n}%
\mathrm{An}(K_{i})$ for any $(\mathcal{D},\mathcal{C})$-subbicomodules $%
K_{1},...,K_{n}\subseteq M$ (i.e. $\mathrm{An}(-)$ in (\ref{An-map-Ke}) and
its restriction in (\ref{An-fi}) are lattice anti-morphisms).

\item $M$ is intrinsically injective.
\end{enumerate}
\end{enumerate}
\end{punto}

\begin{remarks}
\label{duo-duo}Let $M$ be a non-zero $(\mathcal{D},\mathcal{C})$-bicomodule.
If $M$ is self-cogenerator and $^{\mathcal{D}}\mathrm{E}_{M}^{\mathcal{C}}$
is \emph{right-duo} (i.e. every right ideal is a two-sided ideal), then $M$
is duo. On the otherhand, if $M$ is intrinsically injective and $M$ is duo,
then $^{\mathcal{D}}\mathrm{E}_{M}^{\mathcal{C}}$ is right-duo. If $M$ is
self-injective and duo, then every fully invariant $(\mathcal{D},\mathcal{C}%
) $-subbicomodule of $M$ is also duo.
\end{remarks}

\subsection*{Fully coprime (fully cosemiprime) bicomodules}

\begin{punto}
Let $M$ be a non-zero $(\mathcal{D},\mathcal{C})$-bicomodule. For any $R$%
-submodules $X,Y\subseteq M$ we set%
\begin{equation*}
(X\overset{(\mathcal{D},\mathcal{C})}{:_{M}}Y):=\dbigcap \{f^{-1}(Y)\mid
f\in \mathrm{An}_{^{\mathcal{D}}\mathrm{E}_{M}^{\mathcal{C}%
}}(X)\}=\dbigcap\limits_{f\in \mathrm{An}(X)}\{\mathrm{Ker}(\pi _{Y}\circ
f:M\rightarrow M/Y)\}.
\end{equation*}%
If $Y\subseteq M$ is a $(\mathcal{D},\mathcal{C})$-subbicomodule (and $%
f(X)\subseteq X$ for all $f\in $ $^{\mathcal{D}}\mathrm{E}_{M}^{\mathcal{C}}$%
), then $(X\overset{(\mathcal{D},\mathcal{C})}{:_{M}}Y)\subseteq M$ is a
(fully invariant) $(\mathcal{D},\mathcal{C})$-subbicomodule. If $%
X,Y\subseteq M$ are $(\mathcal{D},\mathcal{C})$-subbicomodules, then we call 
$(X\overset{(\mathcal{D},\mathcal{C})}{:_{M}}Y)\subseteq M$ the \emph{%
internal coproduct} of $X$ and $Y$ in $M.$
\end{punto}

\begin{lemma}
\label{inn-ideal}Let $X,Y\subseteq M$ be any $R$-submodules. Then%
\begin{equation}
(X:_{M}^{\mathcal{C}}Y)\subseteq \mathrm{Ke}(\mathrm{An}(X)\circ ^{op}%
\mathrm{An}(Y)),
\end{equation}%
with equality in case $M$ is self-cogenerator and $Y\subseteq M$ is a $(%
\mathcal{D},\mathcal{C})$-subbicomodule.
\end{lemma}

\begin{definition}
Let $M$ be a non-zero $(\mathcal{D},\mathcal{C})$-bicomodule. We call a
non-zero fully invariant $(\mathcal{D},\mathcal{C})$-subbicomodule $0\neq
K\subseteq M:$

\emph{fully }$M$\emph{-coprime}, iff for any fully invariant $(\mathcal{D},%
\mathcal{C})$-subbicomodules $X,Y\subseteq M$ with $K\subseteq (X\overset{(%
\mathcal{D},\mathcal{C})}{:_{M}}Y),$ we have $K\subseteq X$ or $K\subseteq
Y; $

\emph{fully }$M$\emph{-cosemiprime}, iff for any fully invariant $(\mathcal{D%
},\mathcal{C})$-subbicomodule $X\subseteq M$ with $K\subseteq (X\overset{(%
\mathcal{D},\mathcal{C})}{:_{M}}X),$ we have $K\subseteq X;$

In particular, we call $M$ \emph{fully coprime} (\emph{fully cosemiprime}),
iff $M$ is fully $M$-coprime (fully $M$-cosemiprime).
\end{definition}

\subsection*{The fully coprime coradical}

\qquad The \emph{prime spectra }and the associated \emph{prime radicals} for
rings play an important role in the study of structure of rings. Dually, we
define the \emph{fully coprime spectra }and the \emph{fully coprime
coradicals} for bicomodules.

\begin{definition}
Let $M$ be a non-zero $(\mathcal{D},\mathcal{C})$-bicomodule. We define the 
\emph{fully coprime spectrum} of $M$ as%
\begin{equation*}
\mathrm{CPSpec}(M):=\{0\neq K\mid K\subseteq M\text{ is a fully }M\text{%
-coprime }(\mathcal{D},\mathcal{C})\text{-subbicomodule}\}
\end{equation*}%
and the \emph{fully coprime coradical }of $M$ as $\mathrm{CPcorad}%
(M):=\sum\limits_{K\in \mathrm{CPSpec}(M)}K$ ($:=0,$ in case $\mathrm{CPSpec}%
(M)=\varnothing $). Moreover, we set%
\begin{equation*}
\mathrm{CSP}(M):=\{K\mid K\subseteq M\text{ is a fully }M\text{-cosemiprime }%
(\mathcal{D},\mathcal{C})\text{-subbicomodule}\}.
\end{equation*}
\end{definition}

\begin{remark}
We should mention here that the definition of fully coprime (bi)comodules we
present is motivated by the modified version of the definition of coprime
modules (in the sense of Bican et. al. \cite{BJKN80}) as presented in \cite%
{RRW2005}. (Fully) coprime coalgebras over base fields were introduced first
in \cite{NT2001} and considered in \cite{JMR} using the \emph{wedge product}
of subcoalgebras.
\end{remark}

\begin{punto}
Let $M$ be a non-zero $(\mathcal{D},\mathcal{C})$-bicomodule and $L\subseteq
M$ a fully invariant non-zero $(\mathcal{D},\mathcal{C})$-subbicomodule.
Then $L$ is called $\mathrm{E}$\emph{-prime} ($\mathrm{E}$\emph{-semiprime}%
), iff $\mathrm{An}(K)\vartriangleleft $ $^{\mathcal{D}}\mathrm{E}_{M}^{%
\mathcal{C}}$ is prime (semiprime). With $\mathrm{EP}(M)$ ($\mathrm{ESP}(M)$%
) we denote the class of $\mathrm{E}$-prime ($\mathrm{E}$-semiprime) $(%
\mathcal{D},\mathcal{C})$-subbicomodules of $M.$
\end{punto}

The results of \cite{Abu2006} on comodules can be reformulated (with slight
modifications of the proofs) for bicomodules. We state only two of them that
are needed in the sequel.

\begin{proposition}
\label{corad=}Let $M$ be a non-zero $(\mathcal{D},\mathcal{C})$-bicomodule.
If $M$ is self-cogenerator, then $\mathrm{EP}(M)\subseteq \mathrm{CPSpec}(M)$
and $\mathrm{ESP}(M)\subseteq \mathrm{CSP}(M),$ with equality if $M$ is
intrinsically injective. If moreover $^{\mathcal{D}}\mathrm{E}_{M}^{\mathcal{%
C}}$ is right Noetherian, then%
\begin{equation*}
\mathrm{Prad}(^{\mathcal{D}}\mathrm{E}_{M}^{\mathcal{C}})=\mathrm{An}(%
\mathrm{CPcorad}(M))\text{ and }\mathrm{CPcorad}(M)=\mathrm{Ke}(\mathrm{Prad}%
(^{\mathcal{D}}\mathrm{E}_{M}^{\mathcal{C}}));
\end{equation*}%
in particular, $M$ is fully cosemiprime if and only if $M=\mathrm{CPcorad}%
(M).$
\end{proposition}

\begin{proposition}
\label{ro-inn}Let $M$ be a non-zero $(\mathcal{D},\mathcal{C})$-bicomodule
and $0\neq L\subseteq M$ a fully invariant $(\mathcal{D},\mathcal{C})$%
-subbicomodule. If $M$ is self-injective, then%
\begin{equation}
\mathrm{CPSpec}(L)=\mathcal{M}_{f.i.}(L)\cap \mathrm{CPSpec}(M)\text{ and }%
\mathrm{CSP}(L)=\mathcal{M}_{f.i.}(L)\cap \mathrm{CSP}(M);  \label{it<L}
\end{equation}%
hence $\mathrm{CPcorad}(L):=L\cap \mathrm{CPcorad}(M).$
\end{proposition}

\begin{remark}
\label{sub-cop}Let $M$ be a non-zero $(\mathcal{D},\mathcal{C})$-bicomodule.
Then every $L\in \mathcal{S}_{f.i.}(M)$ is trivially a fully coprime $(%
\mathcal{D},\mathcal{C})$-bicomodule. If $M$ is self-injective, then $%
\mathcal{S}_{f.i.}(M)\subseteq \mathrm{CPSpec}(M)$ by Proposition \ref%
{ro-inn}; hence if $M$ has Property $\mathbf{S}_{f.i.},$ then every \emph{%
fully invariant} non-zero $(\mathcal{D},\mathcal{C})$-subbicomodule $%
L\subseteq M$ contains a fully $M$-coprime $(\mathcal{D},\mathcal{C})$%
-subbicomodule $K\subseteq L$ (in particular, $\varnothing \neq \mathrm{%
CPSpec}(L)\subseteq \mathrm{CPSpec}(M)\neq \varnothing $).
\end{remark}

\section{Top Bicomodules}

\qquad In what follows we introduce \emph{top }(\emph{bi})\emph{comodules},
which can be considered (in some sense) as dual to \emph{top }(\emph{bi})%
\emph{modules,} \cite{Lu1999}, \cite{MMS1997}, \cite{Zha1999}. We define a
Zariski topology on the fully coprime spectrum of such (bi)comodules in a
way dual to that of defining the classical Zariski topology on the prime
spectrum of (commutative) rings.

As before, $\mathcal{C}$ is a non-zero $A$-coring and $\mathcal{D}$ is a
non-zero $B$-coring with $_{A}\mathcal{C},\mathcal{\ D}_{B}$ flat. Moreover, 
$M$ is a non-zero $(\mathcal{D},\mathcal{C})$-bicomodule.

\begin{notation}
For every $(\mathcal{D},\mathcal{C})$-subbicomodule $L\subseteq M$ set%
\begin{equation*}
\mathcal{V}_{L}:=\{K\in \mathrm{CPSpec}(M)\mid K\subseteq L\},\text{ }%
\mathcal{X}_{L}:=\{K\in \mathrm{CPSpec}(M)\mid K\nsubseteqq L\}.
\end{equation*}%
Moreover, we set%
\begin{equation*}
\begin{tabular}{lllllll}
$\xi (M)$ & $:=$ & $\{\mathcal{V}_{L}\mid L\in \mathcal{L}(M)\};$ &  & $\xi
_{f.i.}(M)$ & $:=$ & $\{\mathcal{V}_{L}\mid L\in \mathcal{L}_{f.i.}(M)\};$
\\ 
$\tau _{M}$ & $:=$ & $\{\mathcal{X}_{L}\mid L\in \mathcal{L}(M)\};$ &  & $%
\tau _{M}^{f.i.}$ & $:=$ & $\{\mathcal{X}_{L}\mid L\in \mathcal{L}%
_{f.i.}(M)\}.$ \\ 
$\mathbf{Z}_{M}$ & $:=$ & $(\mathrm{CPSpec}(M),\tau _{M});$ &  & $\mathbf{Z}%
_{M}^{f.i.}$ & $:=$ & $(\mathrm{CPSpec}(M),\tau _{M}^{f.i.}).$%
\end{tabular}%
\end{equation*}
\end{notation}

\begin{lemma}
\label{Properties}

\begin{enumerate}
\item $\mathcal{X}_{M}=\emptyset $ and $\mathcal{X}_{\{0\}}=\mathrm{CPSpec}%
(M).$

\item If $\{L_{\lambda }\}_{\Lambda }\subseteq \mathcal{L}(M),$ then $%
\mathcal{X}_{\sum\limits_{\Lambda }L_{\lambda }}\subseteq
\bigcap\limits_{\Lambda }\mathcal{X}_{L_{\lambda }}\subseteq
\bigcup\limits_{\Lambda }\mathcal{X}_{L_{\lambda }}=\mathcal{X}%
_{\bigcap\limits_{\Lambda }L_{\lambda }}.$

\item For any $L_{1},L_{2}\in \mathcal{L}_{f.i.}(M),$ we have $\mathcal{X}%
_{L_{1}+L_{2}}=\mathcal{X}_{L_{1}}\cap \mathcal{X}_{L_{2}}=\mathcal{X}%
_{(L_{1}\overset{(\mathcal{D},\mathcal{C})}{:_{M}}L_{2})}.$
\end{enumerate}
\end{lemma}

\begin{Beweis}
Notice that \textquotedblleft 1\textquotedblright\ and \textquotedblleft
2\textquotedblright\ and the inclusion $\mathcal{X}_{L_{1}+L_{2}}\subseteq 
\mathcal{X}_{L_{1}}\cap \mathcal{X}_{L_{2}}$ in (3) are obvious. If $K\in 
\mathcal{X}_{L_{1}}\cap \mathcal{X}_{L_{2}},$ and $K\notin \mathcal{X}%
_{(L_{1}\overset{(\mathcal{D},\mathcal{C})}{:_{M}}L_{2})},$ then $K\subseteq
L_{1}$ or $K\in L_{2}$ since $K$ is fully $M$-coprime, hence $K\notin 
\mathcal{X}_{L_{1}}$ or $K\notin \mathcal{X}_{L_{2}}$ (a contradiction,
hence $\mathcal{X}_{L_{1}}\cap \mathcal{X}_{L_{2}}\subseteq \mathcal{X}%
_{(L_{1}\overset{(\mathcal{D},\mathcal{C})}{:_{M}}L_{2})}$). Since $%
L_{2}\subseteq M$ is a fully invariant, we have $L_{1}+L_{2}\subseteq (L_{1}%
\overset{(\mathcal{D},\mathcal{C})}{:_{M}}L_{2}),$ hence $\mathcal{X}_{(L_{1}%
\overset{(\mathcal{D},\mathcal{C})}{:_{M}}L_{2})}\subseteq \mathcal{X}%
_{L_{1}+L_{2}}$ and we are done.$\blacksquare $
\end{Beweis}

\begin{remark}
Let $L_{1},L_{2}\subseteq M$ be \emph{arbitrary} $(\mathcal{D},\mathcal{C})$%
-subbicomodules. If $L_{1},L_{2}\subseteq M$ are not fully invariant, then
it is not evident that there exists a $(\mathcal{D},\mathcal{C})$%
-subbicomodule $L\subseteq M$ such that $\mathcal{X}_{L_{1}}\cap \mathcal{X}%
_{L_{2}}=\mathcal{X}_{L}.$ So, for an arbitrary $(\mathcal{D},\mathcal{C})$%
-bicomodule $M,$ the set $\xi (M)$ is not necessarily closed under finite
unions.
\end{remark}

\qquad The remark above motivates the following

\begin{definition}
We call $M$ a \emph{top bicomodule}, iff $\xi (M)$ is closed under finite
unions.
\end{definition}

\qquad As a direct consequence of Lemma \ref{Properties} we get

\begin{theorem}
\label{Topology}$\mathbf{Z}_{M}^{f.i.}:=(\mathrm{CPSpec}(M),\tau
_{M}^{f.i.}) $ is a topological space. In particular, if $M$ is duo, then $M$
is a top $(\mathcal{D},\mathcal{C})$-bicomodule \emph{(}i.e. $\mathbf{Z}%
_{M}:=(\mathrm{CPSpec}(M),\tau _{M})$ is a topological space\emph{)}.
\end{theorem}

\qquad To the end of this section, $M$\ is \emph{duo}, \emph{self-injective}
and $\text{has }$Property $\mathbf{S},\text{ }$so that $\varnothing \neq 
\mathcal{S}(L)=\mathcal{S}_{f.i.}(L)\subseteq \mathrm{CPSpec}(M)$ for every
non-zero $(\mathcal{D},\mathcal{C})$-subbicomodule $0\neq L\subseteq M$ (by
Remark \ref{sub-cop}), and hence a top $(\mathcal{D},\mathcal{C})$%
-bicomodule.

\begin{remarks}
\label{simple-char}Consider the Zariski topology $\mathbf{Z}_{M}:=(\mathrm{%
CPSpec}(M),\tau _{M}).$

\begin{enumerate}
\item $\mathbf{Z}_{M}$ is a $T_{0}$ (Kolmogorov) space.

\item $\mathcal{B}:=\{\mathcal{X}_{L}\mid L\in \mathcal{L}^{f.g.}(M)\}$ is a
basis of open sets for the Zariski topology $\mathbf{Z}_{M}:$ any $K\in 
\mathrm{CPSpec}(M)$ is contained in some $\mathcal{X}_{L}$ for some $L\in 
\mathcal{L}^{f.g.}(M)$ (e.g. $L=0$); and if $L_{1},L_{2}\in \mathcal{L}%
^{f.g.}(M)$ and $K\in \mathcal{X}_{L_{1}}\cap \mathcal{X}_{L_{2}},$ then
setting $L:=L_{1}+L_{2}\in \mathcal{L}^{f.g.}(M),$ we have $K\in \mathcal{X}%
_{L}\subseteq \mathcal{X}_{L_{1}}\cap \mathcal{X}_{L_{2}}.$

\item Let $L\subseteq M$ be a $(\mathcal{D},\mathcal{C}\mathbf{)}$%
-subbicomodule.

\begin{enumerate}
\item $L$ is simple if and only if $L$ is fully $M$-coprime and $\mathcal{V}%
_{L}=\{L\}.$

\item Assume $L\in \mathrm{CPSpec}(M).$ Then $\overline{\{L\}}=\mathcal{V}%
_{L};$ in particular, $L$ is simple if and only if $\{L\}$ is closed in $%
\mathbf{Z}_{M}.$

\item $\mathcal{X}_{L}=\mathrm{CPSpec}(M)$ if and only if $L=0.$

\item If $\mathcal{X}_{L}=\emptyset ,$ then $\mathrm{Corad}(M)\subseteq L.$
\end{enumerate}

\item Let $0\neq L\overset{\theta }{\hookrightarrow }M$ be a non-zero $(%
\mathcal{D},\mathcal{C})$-bicomodule and consider the embedding $\mathrm{%
CPSpec}(L)\overset{\widetilde{\theta }}{\hookrightarrow }\mathrm{CPSpec}(M)$
(compare Proposition \ref{ro-inn}). Since $\theta ^{-1}(\mathcal{V}_{N})=%
\mathcal{V}_{N\cap L}$ for every $N\in \mathcal{L}(M),$ the induced map $%
\mathbf{\theta }:\mathbf{Z}_{L}\rightarrow \mathbf{Z}_{M},$ $K\mapsto \theta
(K)$ is continuous.

\item Let $M\overset{\theta }{\simeq }N$ be an isomorphism of non-zero $(%
\mathcal{D},\mathcal{C})$-bicomodules. Then we have bijections $\mathrm{%
CPSpec}(M)\longleftrightarrow \mathrm{CPSpec}(N)$ and $\mathrm{CSP}%
(M)\longleftrightarrow \mathrm{CSP}(N);$ in particular, $\theta (\mathrm{%
CPcorad}(M))=\mathrm{CPcorad}(N).$ Moreover, $\mathbf{Z}_{M}\approx \mathbf{Z%
}_{N}$ are \emph{homeomorphic spaces}.
\end{enumerate}
\end{remarks}

\begin{theorem}
\label{T1}The following are equivalent:

\begin{enumerate}
\item $\mathrm{CPSpec}(M)=\mathcal{S}(M);$

\item $\mathbf{Z}_{M}$ is discrete;

\item $\mathbf{Z}_{M}$ is a $T_{2}$ \emph{(}Hausdorff space\emph{)} .

\item $\mathbf{Z}_{M}$ is a $T_{1}$ \emph{(}Fr\'{e}cht space\emph{)}.
\end{enumerate}
\end{theorem}

\begin{Beweis}
$(1)\Rightarrow (2).$ For every $K\in \mathrm{CPSpec}(M)=\mathcal{S}(M),$ we
have $\{K\}=\mathcal{X}_{\mathcal{Y}_{K}}$ whence open, where $\mathcal{Y}%
_{K}:=\dsum \{L\in \mathrm{CPSpec}(M)\mid K\nsubseteqq L\}$.

$(2)\Rightarrow (3)$ $\&\ (3)\Rightarrow (4):$ Every discrete topological
space is $T_{2}$ and every $T_{2}$ space is $T_{1}.$

$(4)\Rightarrow (1)$ Let $\mathbf{Z}_{M}$ be $T_{1}$ and suppose $K\in 
\mathrm{CPSpec}(M)\backslash \mathcal{S}(M),$ so that $\{K\}=\mathcal{V}_{L}$
for some $L\in \mathcal{L}(M).$ Since $K$ is not simple, there exists by
assumptions and Remark \ref{sub-cop} $K_{1}\in \mathcal{S}(K)\subseteq 
\mathrm{CPSpec}(M)$ with $K_{1}\subsetneqq K,$ i.e. $\{K_{1},K\}\subsetneqq 
\mathcal{V}_{L}=\{K\},$ a contradiction. Consequently, $\mathrm{CPSpec}(M)=%
\mathcal{S}(M).\blacksquare $
\end{Beweis}

\begin{proposition}
\label{bireg}Let $M$ be self-cogenerator and $^{\mathcal{D}}\mathrm{E}_{M}^{%
\mathcal{C}}$ be Noetherian with every prime ideal maximal \emph{(}e.g. a
biregular ring\footnote{%
a ring in which every two-sided ideal is generated by a central idempotent
(see \cite[3.18(6,7)]{Wis1991}).}\emph{)}.

\begin{enumerate}
\item $\mathcal{S}(M)=\mathrm{CPSpec}(M)$ \emph{(}hence $M$ is subdirectly
irreducible $\Leftrightarrow \left\vert \mathrm{CPSpec}(M)\right\vert =1$%
\emph{)}.

\item If $L\subseteq M$ is a $(\mathcal{D},\mathcal{C})$-subbicomodule, then 
$\mathcal{X}_{L}=\varnothing $ if and only if $\mathrm{Corad}(M)\subseteq L.$
\end{enumerate}
\end{proposition}

\begin{Beweis}
\begin{enumerate}
\item Notice that $\mathcal{S}(M)\subseteq \mathrm{CPSpec}(M)$ by Remark \ref%
{sub-cop}. If $K\in \mathrm{CPSpec}(M),$ then $\mathrm{An}%
(K)\vartriangleleft $ $^{\mathcal{D}}\mathrm{E}_{M}^{\mathcal{C}}$ is prime
by Proposition \ref{corad=}, whence maximal by assumption and it follows
then that $K=\mathrm{Ke}(\mathrm{An}(K))$ is simple (if $0\neq
K_{1}\subsetneqq K,$ for some $K_{1}\in \mathcal{L}(M),$ then $\mathrm{An}%
(K)\subsetneqq \mathrm{An}(K_{1})\subsetneqq $ $^{\mathcal{D}}\mathrm{E}%
_{M}^{\mathcal{C}}$ since $\mathrm{Ke}(-)$ is injective, a contradiction).

\item If $L\subseteq M$ is a $(\mathcal{D},\mathcal{C})$-subbicomodule, then
it follows from \textquotedblleft 1\textquotedblright\ that $\mathcal{X}%
_{L}=\varnothing $ if and only if $\mathrm{Corad}(M)=\mathrm{CPcorad}%
(M)\subseteq L.\blacksquare $
\end{enumerate}
\end{Beweis}

\begin{remark}
\label{PID}Proposition \ref{bireg} corrects \cite[Lemma 2.6.]{NT2001}, which
is absurd since it assumes $C^{\ast }$ PID, while $C$ is not (fully) coprime
(but $C^{\ast }$ domain implies $C$ is (fully) coprime!!).
\end{remark}

\begin{theorem}
\label{compact}If $\left\vert \mathcal{S}(M)\right\vert $ is countable \emph{%
(}finite\emph{), }then $\mathbf{Z}_{M}$ is Lindelof \emph{(}compact\emph{)}.
The converse holds, if $\mathcal{S}(M)=\mathrm{CPSpec}(M).$
\end{theorem}

\begin{Beweis}
Assume $\mathcal{S}(M)=\{S_{\lambda _{k}}\}_{k\geq 1}$ is countable
(finite). Let $\{\mathcal{X}_{L_{\alpha }}\}_{\alpha \in I}$ be an open
cover of $\mathrm{CPSpec}(M)$ (i.e. $\mathrm{CPSpec}(M)\subseteq
\dbigcup\limits_{\alpha \in I}\mathcal{X}_{L_{\alpha }}$). Since $\mathcal{S}%
(M)\subseteq \mathrm{CPSpec}(M)$ we can pick for each $k\geq 1,$ some $%
\alpha _{k}\in I$ such that $S_{\lambda _{k}}\nsubseteqq L_{\alpha _{k}}.$
If $\dbigcap\limits_{k\geq 1}L_{\alpha _{k}}\neq 0,$ then it contains by
Property $\mathbf{S}$ a simple $(\mathcal{D},\mathcal{C})$-subbicomodule $%
0\neq S\subseteq \dbigcap\limits_{k\geq 1}L_{\alpha _{k}},$ (a
contradiction, since $S=S_{\lambda _{k}}\nsubseteqq L_{\alpha _{k}}$ for
some $k\geq 1$). Hence $\dbigcap\limits_{k\geq 1}L_{\alpha _{k}}=0$ and we
conclude that $\mathrm{CPSpec}(M)=\mathcal{X}_{\dbigcap\limits_{k\geq
1}L_{\alpha _{k}}}=\dbigcup\limits_{k\geq 1}\mathcal{X}_{L_{\alpha _{k}}}$
(i.e. $\{\mathcal{X}_{L_{\alpha _{k}}}\mid k\geq 1\}\subseteq \{\mathcal{X}%
_{L_{\alpha }}\}_{\alpha \in I}$ is a countable (finite) subcover). Notice
that if $\mathcal{S}(M)=\mathrm{CPSpec}(M),$ then $\mathbf{Z}_{M}$ is
discrete by Theorem \ref{T1} and so $\mathbf{Z}_{M}$ is Lindelof (compact)
if and only if $\mathrm{CPSpec}(M)$ is countable (finite).$\blacksquare $
\end{Beweis}

\begin{definition}
A collection $\mathcal{G}$ of subsets of a topological space $\mathbf{X}$ is 
\emph{locally finite}, iff every point of $\mathbf{X}$ has a neighbourhood
that intersects only finitely many elements of $\mathcal{G}.$
\end{definition}

\begin{proposition}
\label{lf}Let $\mathcal{K}=\{K_{\lambda }\}_{\Lambda }\subseteq \mathcal{S}%
(M)$ be a non-empty family of simple $(\mathcal{D},\mathcal{C})$%
-subbicomodules. If $\left\vert \mathcal{S}(L)\right\vert <\infty $ for
every $L\in \mathrm{CPSpec}(M),$ then $\mathcal{K}$ is locally finite.
\end{proposition}

\begin{Beweis}
Let $L\in \mathrm{CPSpec}(M)$ and set $F:=\dsum \{K\in \mathcal{K}\mid
K\nsubseteqq L\}.$ Since $\left\vert \mathcal{S}(L)\right\vert <\infty ,$
there exists a finite number of simple $(\mathcal{D},\mathcal{C})$%
-subbicomodules $\{S_{\lambda _{1}},..,S_{\lambda _{n}}\}=\mathcal{K}\cap 
\mathcal{V}_{L}.$ If $L\subseteq F,$ then $0\neq L\subseteq
\dsum\limits_{i=1}^{n}S_{\lambda _{i}}\subseteq (S_{\lambda _{1}}:_{M}^{(%
\mathcal{D},\mathcal{C})}\dsum\limits_{i=2}^{n}S_{\lambda _{i}})$ and it
follows by induction that $0\neq L\subsetneqq S_{\lambda _{i}}$ for some $%
1\leq i\leq n$ (a contradiction, since $S_{\lambda _{i}}$ is simple), whence 
$L\in \mathcal{X}_{F}.$ It is clear then that $\mathcal{K}\cap \mathcal{X}%
_{F}=\{K_{\lambda _{1}},..,K_{\lambda _{n}}\}$ and we are done.$\blacksquare 
$
\end{Beweis}

\begin{definition}
(\cite{Bou1966}, \cite{Bou1998})\ A topological space $\mathbf{X}$ is said
to be \emph{irreducible} (\emph{connected}), iff $\mathbf{X}$ is not the
(disjoint) union of two proper closed subsets; equivalently, iff the
intersection of any two non-empty open subsets is non-empty (the only
subsets of $\mathbf{X}$ that are open and closed are $\varnothing $ and $%
\mathbf{X}$). A maximal irreducible subspace of $\mathbf{X}$ is called an 
\emph{irreducible component}.
\end{definition}

\begin{proposition}
\label{duo-irr}$\mathrm{CPSpec}(M)$ is irreducible if and only if $\mathrm{%
CPcorad}(M)$ is fully $M$-coprime.
\end{proposition}

\begin{Beweis}
Let $\mathrm{CPSpec}(M)$ be irreducible. By Remark \ref{sub-cop}, $\mathrm{%
CPcorad}(M)\neq 0.$ Suppose that $\mathrm{CPcorad}(M)$ is not fully $M$%
-coprime, so that there exist $(\mathcal{D},\mathcal{C})$-subbicomodules $%
X,Y\subseteq M$ with $\mathrm{CPcorad}(M)\subseteq (X:_{M}^{(\mathcal{D},%
\mathcal{C})}Y)$ but $\mathrm{CPcorad}(M)\nsubseteqq X$ and $\mathrm{CPcorad}%
(M)\nsubseteqq Y.$ It follows then that $\mathrm{CPSpec}(M)=\mathcal{V}%
_{(X:_{M}^{(\mathcal{D},\mathcal{C})}Y)}=\mathcal{V}_{X}\cup \mathcal{V}_{Y}$
a union of proper closed subsets, a contradiction. Consequently, $\mathrm{%
CPcorad}(M)$ is fully $M$-coprime.

On the otherhand, assume $\mathrm{CPcorad}(M)\in \mathrm{CPSpec}(M)$ and
suppose that $\mathrm{CPSpec}(M)=\mathcal{V}_{L_{1}}\cup \mathcal{V}_{L_{2}}=%
\mathcal{V}_{(L_{1}:_{M}^{(\mathcal{D},\mathcal{C})}L_{2})}$ for some $(%
\mathcal{D},\mathcal{C})$-subbicomodules $L_{1},L_{2}\subseteq M.$ It
follows then that $\mathrm{CPcorad}(M)\subseteq L_{1},$ so that $\mathcal{V}%
_{L_{1}}=\mathrm{CPSpec}(M);$ or $\mathrm{CPcorad}(M)\subseteq L_{2},$ so
that $\mathcal{V}_{L_{2}}=\mathrm{CPSpec}(M).$ Consequently $\mathrm{CPSpec}%
(M)$ is not the union of two \emph{proper} closed subsets, i.e. it is
irreducible.$\blacksquare $
\end{Beweis}

\begin{lemma}
\label{it-irr}

\begin{enumerate}
\item $M$ is subdirectly irreducible if and only if the intersection of any
two non-empty closed subsets of $\mathrm{CPSpec}(M)$ is non-empty.

\item If $M$ is subdirectly irreducible, then $\mathrm{CPSpec}(M)$ is
connected. If $\mathrm{CPSpec}(M)$ is connected and $\mathrm{CPSpec}(M)=%
\mathcal{S}(M),$ then $M$ is subdirectly irreducible.
\end{enumerate}
\end{lemma}

\begin{Beweis}
\begin{enumerate}
\item Assume that $M$ is subdirectly irreducible with unique simple $(%
\mathcal{D},\mathcal{C})$-subbicomodule $0\neq S\subseteq M.$ If $\mathcal{V}%
_{L_{1}},$ $\mathcal{V}_{L_{2}}\subseteq \mathrm{CPSpec}(M)$ are any two
non-empty closed subsets, then $L_{1}\neq 0\neq L_{2}$ and so $\mathcal{V}%
_{L_{1}}\cap \mathcal{V}_{L_{2}}=\mathcal{V}_{L_{1}\cap L_{2}}\neq
\varnothing ,$ since $S\subseteq L_{1}\cap L_{2}\neq 0.$ On the otherhand,
assume that the intersection of any two non-empty closed subsets of $\mathrm{%
CPSpec}(M)$ is non-empty. Let $0\neq L_{1},L_{2}\subseteq M$ be any non-zero 
$(\mathcal{D},\mathcal{C})$-subbicomodules, so that $\mathcal{V}_{L_{1}}\neq
\varnothing \neq \mathcal{V}_{L_{2}}.$ By assumption $\mathcal{V}_{L_{1}\cap
L_{2}}=\mathcal{V}_{L_{1}}\cap \mathcal{V}_{L_{2}}\neq \varnothing ,$ hence $%
L_{1}\cap L_{2}\neq 0$ and it follows by \ref{ft=2} that $M$ is subdirectly
irreducible.

\item If $M$ is subdirectly irreducible, then $\mathrm{CPSpec}(M)$ is
connected by \textquotedblleft 1\textquotedblright . On the otherhand, if $%
\mathrm{CPSpec}(M)=\mathcal{S}(M),$ then $\mathbf{Z}_{M}$ is discrete by
Theorem \ref{T1} and so $M$ is subdirectly irreducible (since a discrete
topological space is connected if and only if it has only one point).$%
\blacksquare $
\end{enumerate}
\end{Beweis}

\begin{proposition}
\begin{enumerate}
\item If $K\in \mathrm{CPSpec}(M),$ then $\mathcal{V}_{K}\subseteq \mathrm{%
CPSpec}(M)$ is irreducible.

\item If $\mathcal{V}_{L}$ is an irreducible component of $\mathbf{Z}_{M},$
then $L$ is a maximal fully $M$-coprime $(\mathcal{D},\mathcal{C})$%
-subbicomodule.
\end{enumerate}
\end{proposition}

\begin{Beweis}
\begin{enumerate}
\item Let $K\in \mathrm{CPSpec}(M)$ and suppose $\mathcal{V}_{K}=A\cup B=(%
\mathcal{V}_{K}\cap \mathcal{V}_{X})\cup (\mathcal{V}_{K}\cap \mathcal{V}%
_{Y})$ for two $(\mathcal{D},\mathcal{C})$-subbicomodules $X,Y\subseteq M$
(so that $A,B\subseteq \mathcal{V}_{K}$ are closed subsets w.r.t. the
relative topology on $\mathcal{V}_{K}\hookrightarrow \mathrm{CPSpec}(M)$).
It follows then that $\mathcal{V}_{K}=(\mathcal{V}_{K\cap X})\cup (\mathcal{V%
}_{K\cap Y})=\mathcal{V}_{(K\cap X:_{M}^{(\mathcal{D},\mathcal{C})}K\cap Y)}$
and so $K\subseteq (K\cap X:_{M}^{(\mathcal{D},\mathcal{C})}K\cap Y),$ hence 
$K\subseteq X$ so that $\mathcal{V}_{K}=A;$ or $K\subseteq Y,$ so that $%
\mathcal{V}_{K}=B.$ Consequently $\mathcal{V}_{K}$ is irreducible.

\item Assume $\mathcal{V}_{L}$ is an irreducible component of $\mathrm{CPSpec%
}(M)$ for some $0\neq L\in \mathcal{L}(M).$ If $L\subseteq K$ for some $K\in 
\mathrm{CPSpec}(M),$ then $\mathcal{V}_{L}\subseteq \mathcal{V}_{K}$ and it
follows then that $L=K$ (since $\mathcal{V}_{K}\subseteq \mathrm{CPSpec}(M)$
is irreducible by \textquotedblleft 1\textquotedblright ). We conclude then
that $L$ is fully $M$-coprime and is moreover maximal in $\mathrm{CPSpec}%
(M).\blacksquare $
\end{enumerate}
\end{Beweis}

\begin{lemma}
\label{1n}If $n\geq 2$ and $\mathcal{A}=\{K_{1},...,K_{n}\}\subseteq \mathrm{%
CPSpec}(M)$ is a connected subset, then for every $i\in \{1,...,n\},$ there
exists $j\in \{1,...,n\}\backslash \{i\}$ such that $K_{i}\subseteq K_{j}$
or $K_{j}\subseteq K_{i}.$
\end{lemma}

\begin{Beweis}
Without loss of generality, suppose $K_{1}\nsubseteqq K_{j}$ and $%
K_{j}\nsubseteqq K_{1}$ for all $2\leq j\leq n$ and set $F:=\dsum%
\limits_{i=2}^{n}K_{i},$ $W_{1}:=\mathcal{A}\cap \mathcal{X}%
_{K_{1}}=\{K_{2},...,K_{n}$ and $W_{2}:=\mathcal{A}\cap \mathcal{X}%
_{F}=\{K_{1}\}$ (if $n=2,$ then clearly $W_{2}=\{K_{1}\};$ if $n>2$ and $%
K_{1}\notin W_{2},$ then $K_{1}\subseteq
\dsum\limits_{i=2}^{n}K_{i}\subseteq (K_{2}:_{M}^{(\mathcal{D},\mathcal{C}%
)}\dsum\limits_{i=3}^{n}K_{i})$ and it follows that $K_{1}\subseteq
\dsum\limits_{i=3}^{n}K_{i};$ by induction one shows that $K_{1}\subseteq
K_{n},$ a contradiction). So $\mathcal{A}=W_{1}\cup W_{2},$ a disjoint union
of proper non-empty open subsets (a contradiction).$\blacksquare $
\end{Beweis}

\begin{notation}
For $\mathcal{A}\subseteq \mathrm{CPSpec}(M)$ set $\varphi (\mathcal{A}%
):=\dsum_{K\in \mathcal{A}}K$ ($:=0,$ iff $\mathcal{A}=\varnothing $).
Moreover, set%
\begin{equation*}
\mathbf{CL}(\mathbf{Z}_{M}):=\{\mathcal{A}\subseteq \mathrm{CPSpec}(M)\mid 
\mathcal{A}=\overline{\mathcal{A}}\}\text{ and }\mathcal{E}(M):=\{L\in 
\mathcal{L}(M)\mid \mathrm{CPcorad}(L)=L\}.
\end{equation*}
\end{notation}

\begin{lemma}
The closure of any subset $\mathcal{A}\subseteq \mathrm{CPSpec}(M)$ is $%
\overline{\mathcal{A}}=\mathcal{V}_{\varphi (\mathcal{A})}.$
\end{lemma}

\begin{Beweis}
Let $\mathcal{A}\subseteq \mathrm{CPSpec}(M).$ Since $\mathcal{A}\subseteq 
\mathcal{V}_{\varphi (\mathcal{A})}$ and $\mathcal{V}_{\varphi (\mathcal{A}%
)} $ is a closed set, we have $\overline{\mathcal{A}}\subseteq \mathcal{V}%
_{\varphi (\mathcal{A})}.$ On the other hand, suppose $H\in \mathcal{V}%
_{\varphi (\mathcal{A})}\backslash \mathcal{A}$ and let $\mathcal{X}_{L}$ be
a neighbourhood of $H,$ so that $H\nsubseteqq L.$ Then there exists $W\in 
\mathcal{A}$ with $W\nsubseteqq L$ (otherwise $H\subseteq \varphi (\mathcal{A%
})\subseteq L,$ a contradiction), i.e. $W\in \mathcal{X}_{L}\cap (\mathcal{A}%
\backslash \{H\})\neq \varnothing $ and so $K$ is a cluster point of $%
\mathcal{A}\mathfrak{.}$ Consequently, $\overline{\mathcal{A}}=\mathcal{V}%
_{\varphi (\mathcal{A})}.\blacksquare $
\end{Beweis}

\begin{theorem}
\label{11}We have a bijection $\mathbf{CL}(\mathbf{Z}_{M})%
\longleftrightarrow \mathcal{E}(M).$ If $M$ is self-cogenerator and $^{%
\mathcal{D}}\mathrm{E}_{M}^{\mathcal{C}}$ is right Noetherian, then there is
a bijection $\mathbf{CL}(\mathbf{Z}_{M})\backslash \{\varnothing
\}\longleftrightarrow \mathrm{CSPSpec}(M).$
\end{theorem}

\begin{Beweis}
For $L\in \mathcal{E}(M),$ set $\psi (L):=\mathcal{V}_{L}.$ Then for $L\in 
\mathcal{E}(M)$ and $\mathcal{A}\in \mathbf{CL}(\mathbf{Z}_{M})$ we have $%
\varphi (\psi (L))=\varphi (\mathcal{V}_{L})=L\cap \mathrm{C}\mathrm{Pcorad}%
(M)=\mathrm{C}\mathrm{Pcorad}(L)=L$ and $\psi (\varphi (\mathcal{A}))=%
\mathcal{V}_{\varphi (\mathcal{A})}=\overline{\mathcal{A}}=\mathcal{A}.$ If $%
M$ is self-cogenerator and $^{\mathcal{D}}\mathrm{E}_{M}^{\mathcal{C}}$ is
right Noetherian, then $\mathrm{CSPSpec}(M)=\mathcal{E}(M)\backslash \{0\}$
by Proposition \ref{corad=} and we are done.$\blacksquare $
\end{Beweis}

\section{Applications and Examples}

\qquad In this section we give some applications and examples. First of all
we remark that taking $\mathcal{D}:=R$ ($\mathcal{C}:=R$), considered with
the trivial coring structure, our results on the Zariski topology for
bicomodules in the third section can be reformulated for Zariski topology on
the fully coprime spectrum of right $\mathcal{C}$-comodules (left $\mathcal{D%
}$-comodules). However, our main application will be to the Zariski topology
on the fully coprime spectrum of non-zero corings, considered as \emph{duo
bicomodules} in the canonical way.

Throughout this section, $\mathcal{C}$ is a non-zero $A$-coring with $_{A}%
\mathcal{C}$ and $\mathcal{C}_{A}$ flat.

\begin{punto}
The $(A,A)$-bimodule $^{\ast }\mathcal{C}^{\ast }:=\mathrm{Hom}_{(A,A)}(%
\mathcal{C},A):=$ $^{\ast }\mathcal{C}\cap \mathcal{C}^{\ast }$ is an $%
A^{op} $-ring with multiplication $(f\ast g)(c)=\dsum f(c_{1})g(c_{2})$ for
all $f,g\in $ $^{\ast }\mathcal{C}^{\ast }$ and unit $\varepsilon _{\mathcal{%
C}};$ hence every $(\mathcal{C},\mathcal{C})$-bicomodule $M$ is a $(^{\ast }%
\mathcal{C}^{\ast },^{\ast }\mathcal{C}^{\ast })$-bimodule and the \emph{%
centralizer }%
\begin{equation*}
\mathbf{C}(M):=\{f\in \text{ }^{\ast }\mathcal{C}^{\ast }\mid
f\rightharpoonup m=m\leftharpoonup f\text{ for all }m\in M\}
\end{equation*}%
is an $R$-algebra. If $M$ is faithful as a left (right) $^{\ast }\mathcal{C}%
^{\ast }$-module, then $\mathbf{C}(M)\subseteq Z(^{\ast }\mathcal{C}^{\ast
}).$
\end{punto}

\begin{punto}
Considering $\mathcal{C}$ as a $(\mathcal{C},\mathcal{C})$-bicomodule in the
natural way, $\mathcal{C}$ is a $(^{\ast }\mathcal{C}^{\ast },^{\ast }%
\mathcal{C}^{\ast })$-bimodule that is faithful as a left (right) $^{\ast }%
\mathcal{C}^{\ast }$-module, hence the centralizer%
\begin{equation*}
\mathbf{C}(\mathcal{C}):=\{f\in \text{ }^{\ast }\mathcal{C}^{\ast }\mid
f\rightharpoonup c=c\leftharpoonup f\text{ for every }c\in \mathcal{C}\}
\end{equation*}%
embeds in the center of $^{\ast }\mathcal{C}^{\ast }$ as an $R$-subalgebra,
i.e. $\mathbf{C}(\mathcal{C})\hookrightarrow Z(^{\ast }\mathcal{C}^{\ast }).$
If $ac=ca$ for all $a\in A,$ then we have a morphism of $R$-algebras $\eta
:Z(A)\rightarrow \mathbf{C}(\mathcal{C}),$ $a\mapsto \lbrack \varepsilon _{%
\mathcal{C}}(a-)=\varepsilon _{\mathcal{C}}(-a)].$
\end{punto}

\begin{remark}
Notice that $\mathbf{C}(\mathcal{C})\subseteq Z(^{\ast }\mathcal{C}%
)\subseteq Z(^{\ast }\mathcal{C}^{\ast })$ and $\mathbf{C}(\mathcal{C}%
)\subseteq Z(\mathcal{C}^{\ast })\subseteq Z(^{\ast }\mathcal{C}^{\ast })$
(compare \cite[17.8. (4)]{BW2003}). If $_{A}\mathcal{C}$ ($\mathcal{C}_{A}$)
is $A$-cogenerated, then it follows by \cite[19.10 (3)]{BW2003} that $%
Z(^{\ast }\mathcal{C})=\mathbf{C}(\mathcal{C})\subseteq Z(\mathcal{C}^{\ast
})$ ($Z(\mathcal{C}^{\ast })=\mathbf{C}(\mathcal{C})\subseteq Z(^{\ast }%
\mathcal{C})$). If $_{A}\mathcal{C}_{A}$ is $A$-cogenerated, then $Z(^{\ast }%
\mathcal{C}^{\ast })\subseteq \mathbf{C}(\mathcal{C})$ (e.g. \cite[19.10 (4)]%
{BW2003}), whence $Z(^{\ast }\mathcal{C})=Z(^{\ast }\mathcal{C}^{\ast })=Z(%
\mathcal{C}^{\ast }).$
\end{remark}

\begin{lemma}
\label{phi-M}For every $(\mathcal{C},\mathcal{C})$-bicomodule $M$ we have a
morphism of $R$-algebras%
\begin{equation}
\phi _{M}:\mathbf{C}(M)\rightarrow \text{ }^{\mathcal{C}}\mathrm{End}^{%
\mathcal{C}}(M)^{op},\text{ }f\mapsto \lbrack m\mapsto f\rightharpoonup
m=m\leftharpoonup f]\text{ \emph{(}with }\mathrm{Im}(\phi _{M})\subseteq Z(^{%
\mathcal{C}}\mathrm{E}_{M}^{\mathcal{C}})\text{\emph{)}.}  \label{E=Z}
\end{equation}
\end{lemma}

\begin{Beweis}
First of all we prove that $\phi _{M}$ is well-defined: for $f\in \mathbf{C}%
(M)$ and $m\in M$ we have%
\begin{equation*}
\begin{tabular}{lll}
$\sum (\phi _{M}(f)(m))_{<0>}\otimes _{A}(\phi _{M}(f)(m))_{<1>}$ & $=$ & $%
\sum (m\leftharpoonup f)_{<0>}\otimes _{A}(m\leftharpoonup f)_{<1>}$ \\ 
& $=$ & $\sum f(m_{<-1>})m_{<0><0>}\otimes _{A}m_{<0><1>}$ \\ 
& $=$ & $\sum f(m_{<0><-1>})m_{<0><0>}\otimes _{A}m_{<1>}$ \\ 
& $=$ & $\sum (m_{<0>}\leftharpoonup f)\otimes _{A}m_{<1>}$ \\ 
& $=$ & $\sum \phi _{M}(f)(m_{<0>})\otimes _{A}m_{<1>},$%
\end{tabular}%
\end{equation*}%
and%
\begin{equation*}
\begin{tabular}{lll}
$\sum (\phi _{M}(f)(m))_{<-1>}\otimes _{A}(\phi _{M}(f)(m))_{<0>}$ & $=$ & $%
\sum (f\rightharpoonup m)_{<-1>}\otimes _{A}(f\rightharpoonup m)_{<0>}$ \\ 
& $=$ & $\sum m_{<0><-1>}\otimes _{A}m_{<0><0>}f(m_{<1>})$ \\ 
& $=$ & $\sum m_{<-1>}\otimes _{A}m_{<0><0>}f(m_{<0><1>})$ \\ 
& $=$ & $\sum m_{<-1>}\otimes _{A}(f\rightharpoonup m_{<0>})$ \\ 
& $=$ & $\sum m_{<-1>}\otimes _{A}\phi _{M}(f)(m_{<0>}),$%
\end{tabular}%
\end{equation*}%
i.e. $\phi _{M}(f):M\rightarrow M$ is $(\mathcal{C},\mathcal{C})$%
-bicolinear. Obviously, $\phi _{M}(f\ast g)=\phi _{M}(f)\circ ^{op}\phi
_{M}(g)$ for all $f,g\in \mathbf{C}(M),$ i.e. $\phi _{M}$ is a morphism of $%
R $-algebras. Moreover, since every $g\in $ $^{\mathcal{C}}\mathrm{E}_{M}^{%
\mathcal{C}}$ is $(^{\ast }\mathcal{C}^{\ast },^{\ast }\mathcal{C}^{\ast })$%
-bilinear, we have $g(f\rightharpoonup m)=f\rightharpoonup g(m)$ for every $%
f\in $ $^{\ast }\mathcal{C}^{\ast }$ and $m\in M,$ i.e. $\mathrm{Im}(\phi
_{M})\subseteq Z(^{\mathcal{C}}\mathrm{E}_{M}^{\mathcal{C}}).\blacksquare $
\end{Beweis}

\begin{lemma}
\label{End(C)=Z}We have an isomorphism of $R$-algebras $\mathbf{C}(\mathcal{C%
})\overset{\phi _{\mathcal{C}}}{\simeq }$ $^{\mathcal{C}}\mathrm{End}^{%
\mathcal{C}}(\mathcal{C}),$ with inverse $\psi _{\mathcal{C}}:g\mapsto
\varepsilon _{\mathcal{C}}\circ g.$ In particular, $(^{\mathcal{C}}\mathrm{%
End}^{\mathcal{C}}(\mathcal{C}),\circ )$ is commutative and $\mathcal{C}\in $
$^{\mathcal{C}}\mathbb{M}^{\mathcal{C}}$ is duo.
\end{lemma}

\begin{Beweis}
First of all we prove that $\psi $ is well-defined: for $g\in $ $^{\mathcal{C%
}}\mathrm{End}^{\mathcal{C}}(\mathcal{C})$ and $c\in \mathcal{C}$ we have%
\begin{equation*}
\begin{tabular}{lllllll}
$\psi _{\mathcal{C}}(g)\rightharpoonup c$ & $=$ & $\sum c_{1}\psi (g)(c_{2})$
& $=$ & $\sum c_{1}\varepsilon _{\mathcal{C}}(g(c_{2}))$ & $=$ & $\sum
g(c)_{1}\varepsilon _{\mathcal{C}}(g(c)_{2})$ \\ 
& $=$ & $g(c)$ & $=$ & $\sum \varepsilon _{\mathcal{C}}(g(c)_{1})g(c)_{2}$ & 
$=$ & $\sum \varepsilon _{\mathcal{C}}(g(c_{1}))c_{2}$ \\ 
& $=$ & $\sum \psi (g)(c_{1})c_{2}$ & $=$ & $c\leftharpoonup \psi (g),$ &  & 
\end{tabular}%
\end{equation*}%
i.e. $\psi _{\mathcal{C}}(g)\in \mathbf{C}(\mathcal{C}).$ For any $f\in 
\mathbf{C}(\mathcal{C}),$ $g\in $ $^{\mathcal{C}}\mathrm{End}^{\mathcal{C}}(%
\mathcal{C})$ and $c\in \mathcal{C}$ we have $((\psi _{\mathcal{C}}\circ
\phi _{\mathcal{C}})(f))(c)=\varepsilon _{\mathcal{C}}(\phi _{\mathcal{C}%
}(f)(c))=\varepsilon _{\mathcal{C}}(f\rightharpoonup c)=f(c)$ and $((\phi _{%
\mathcal{C}}\circ \psi _{\mathcal{C}})(g))(c)=\sum c_{1}\psi _{\mathcal{C}%
}(g)(c_{2})=\sum c_{1}\varepsilon _{\mathcal{C}}(g(c_{2}))=\sum
g(c)_{1}\varepsilon _{\mathcal{C}}(g(c)_{2})=g(c).\blacksquare $
\end{Beweis}

\subsection*{Zariski topologies for corings}

\begin{definition}
A right (left) $\mathcal{C}$-subcomodule $K\subseteq \mathcal{C}$ is called
a \emph{right }(\emph{left})\emph{\ }$\mathcal{C}$\emph{-coideal}. A $(%
\mathcal{C},\mathcal{C})$-subbicomodule of $M$ is called a $\mathcal{C}$%
\emph{-bicoideal}.
\end{definition}

\begin{notation}
With $\mathcal{B}(\mathcal{C})$ we denote the class of $\mathcal{C}$%
-bicoideals and with $\mathcal{L}(\mathcal{C}^{r})$ (resp. $\mathcal{L}(%
\mathcal{C}^{l})$) the class of right (left) $\mathcal{C}$-coideals. For a $%
\mathcal{C}$-bicoideal $K\in \mathcal{B}(\mathcal{C}),$ $K^{r}$ ($K^{l}$)
indicates that we consider $K$ as a right (left) $\mathcal{C}$-comodule,
rather than a $(\mathcal{C},\mathcal{C})$-bicomodule. We also set%
\begin{equation*}
\begin{tabular}{llllll}
$\mathrm{CPSpec}(\mathcal{C})$ & $:=$ & $\{K\in \mathcal{B}(\mathcal{C})\mid
K$ is fully $\mathcal{C}$-coprime$\};$ & $\tau _{\mathcal{C}}$ & $:=$ & $\{%
\mathcal{X}_{L}\mid L\in \mathcal{B}(\mathcal{C})\};$ \\ 
$\mathrm{CPSpec}(\mathcal{C}^{r})$ & $:=$ & $\{K\in \mathcal{B}(\mathcal{C}%
)\mid K^{r}$ is fully $\mathcal{C}^{r}$-coprime$\};$ & $\tau _{\mathcal{C}%
^{r}}$ & $:=$ & $\{\mathcal{X}_{L}\mid L\in \mathcal{L}(\mathcal{C}^{r})\};$
\\ 
$\mathrm{CPSpec}(\mathcal{C}^{l})$ & $:=$ & $\{K\in \mathcal{B}(\mathcal{C}%
)\mid K^{l}$ is fully $\mathcal{C}^{l}$-coprime$\};$ & $\tau _{\mathcal{C}%
^{l}}$ & $:=$ & $\{\mathcal{X}_{L}\mid L\in \mathcal{L}(\mathcal{C}^{l})\}.$%
\end{tabular}%
\end{equation*}
\end{notation}

In what follows we announce only the main result on the Zariski topologies
for corings, leaving to the interested reader the restatement of the other
results of the third section.

\begin{theorem}
\begin{enumerate}
\item $\mathbf{Z}_{\mathcal{C}}:=(\mathrm{CPSpec}(\mathcal{C}),\tau _{%
\mathcal{C}})$ is a topological space.

\item $\mathbf{Z}_{\mathcal{C}r}^{f.i.}:=(\mathrm{CPSpec}(C^{r}),\tau _{%
\mathcal{C}^{r}}^{f.i.})$ and $\mathbf{Z}_{\mathcal{C}^{l}}^{f.i.}:=(\mathrm{%
CPSpec}(C^{l}),\tau _{\mathcal{C}^{l}}^{f.i.})$ are topological spaces.
\end{enumerate}
\end{theorem}

\begin{proposition}
\label{th-tel}Let $\theta :\mathcal{C}\rightarrow \mathcal{C}^{\prime }$ be
a morphism of non-zero $A$-corings with $_{A}\mathcal{C},$ $_{A}\mathcal{C}%
^{\prime }$ flat, $\mathcal{C}^{r}$ intrinsically injective self-cogenerator
and $\mathcal{C}^{\prime r}$ self-cogenerator.

\begin{enumerate}
\item If $\theta $ is injective and $\mathcal{C}^{\prime r}$ is
self-injective, or if $\mathcal{C}^{\ast }$ is right-duo, then we have a map 
$\widetilde{\theta }:\mathrm{CPSpec}(\mathcal{C}^{r})\rightarrow \mathrm{%
CPSpec}(\mathcal{C}^{\prime r}),$ $K\mapsto \theta (K)$ \emph{(}and so $%
\theta (\mathrm{CPcorad}(\mathcal{C}^{r}))\subseteq \mathrm{CPcorad}(%
\mathcal{C}^{\prime r})$\emph{)}.

\item If $\mathcal{C}^{r},$ $\mathcal{C}^{\prime r}$ are duo, then the
induced map $\mathbf{\theta }:\mathbf{Z}_{\mathcal{C}^{r}}\rightarrow 
\mathbf{Z}_{\mathcal{C}^{\prime r}}$ is continuous.

\item If every $K\in \mathrm{CPSpec}(\mathcal{C}^{r})$ is inverse image of a 
$K^{\prime }\in \mathrm{CPSpec}(\mathcal{C}^{\prime r}),$ then $\widetilde{%
\theta }$ is injective.

\item If $\theta $ is injective and $\mathcal{C}^{\prime r}$ is
self-injective, then $\mathbf{\theta }:\mathbf{Z}_{\mathcal{C}%
^{r}}^{f.i.}\rightarrow \mathbf{Z}_{\mathcal{C}^{\prime r}}^{f.i.}$ is
continuous. If moreover, $\widetilde{\theta }:\mathrm{CPSpec}(\mathcal{C}%
^{r})\rightarrow \mathrm{CPSpec}(\mathcal{C}^{\prime r})$ is surjective,
then $\mathbf{\theta }$ is open and closed.

\item If $\mathcal{C}\overset{\theta }{\simeq }\mathcal{C}^{\prime },$ then $%
\mathbf{Z}_{\mathcal{C}^{r}}^{f.i.}\overset{\mathbf{\theta }}{\approx }%
\mathbf{Z}_{\mathcal{C}^{\prime r}}^{f.i.}$ \emph{(}homeomorphic spaces\emph{%
)}.
\end{enumerate}
\end{proposition}

\begin{Beweis}
First of all notice for every $K\in \mathcal{L}(\mathcal{C}^{r})$ ($K\in 
\mathcal{B}(\mathcal{C})$), we have $\theta (K)\in \mathcal{L}(\mathcal{C}%
^{\prime r})$ ($\theta (K)\in \mathcal{B}(\mathcal{C}^{\prime })$) and for
every $K^{\prime }\in \mathcal{L}(\mathcal{C}^{\prime r})$ ($K^{\prime }\in 
\mathcal{B}(\mathcal{C}^{\prime })$), $\theta ^{-1}(K^{\prime })\in \mathcal{%
L}(\mathcal{C}^{r})$ ($\theta ^{-1}(K^{\prime })\in \mathcal{B}(\mathcal{C})$%
).

\begin{enumerate}
\item If $\theta $ is injective and $\mathcal{C}^{\prime r}$ is
self-injective, then $\mathrm{CPSpec}(\mathcal{C}^{r})=\mathcal{B}(\mathcal{C%
})\cap \mathrm{C}\mathrm{PSpec}(\mathcal{C}^{\prime r})$ by \cite[%
Proposition 4.7.]{Abu2006}. Assume now that $\mathcal{C}^{\ast }$ is
right-duo. Since $\theta $ is a morphism of $A$-corings, the canonical map $%
\theta ^{\ast }:\mathcal{C}^{\prime \ast }\rightarrow \mathcal{C}^{\ast }$
is a morphism of $A^{op}$-rings. If $K\in \mathrm{CPSpec}(\mathcal{C}^{r}),$
then $\mathrm{ann}_{\mathcal{C}^{\ast }}(K)\vartriangleleft \mathcal{C}%
^{\ast }$ is a prime ideal by \cite[Proposition 4.10.]{Abu2006}, whence
completely prime since $\mathcal{C}^{\ast }$ is right-duo. It follows then
that $\mathrm{ann}_{\mathcal{C}^{\prime \ast }}(\theta (K))=\theta (K)^{\bot 
\mathcal{C}^{\prime \ast }}=(\theta ^{\ast })^{-1}(K^{\bot \mathcal{C}^{\ast
}})=(\theta ^{\ast })^{-1}(\mathrm{ann}_{\mathcal{C}^{\ast }}(K))$ is a
prime ideal, whence $\theta (K)\in \mathrm{CPSpec}(\mathcal{C}^{\prime r})$
by \cite[Proposition 4.10.]{Abu2006}. It is obvious then that $\theta (%
\mathrm{CPcorad}(\mathcal{C}^{r}))\subseteq \mathrm{CPcorad}(\mathcal{C}%
^{\prime r}).$

\item Since $\mathcal{C}^{r}\in \mathbb{M}^{\mathcal{C}},$ $\mathcal{C}%
^{\prime r}\in \mathbb{M}^{\mathcal{C}^{\prime }}$ are duo, $\mathbf{Z}_{%
\mathcal{C}^{r}}:=\mathbf{Z}_{\mathcal{C}^{r}}^{f.i.}$ and $\mathbf{Z}_{%
\mathcal{C}^{\prime r}}:=\mathbf{Z}_{\mathcal{C}^{\prime r}}^{f.i.}$ are
topological spaces. Since $\mathcal{C}^{r}$ is intrinsically injective, $%
\mathcal{C}^{\ast }$ is right-due and by \textquotedblleft
1\textquotedblright\ $\widetilde{\theta }:\mathrm{CPSpec}(\mathcal{C}%
^{r})\rightarrow \mathrm{CPSpec}(\mathcal{C}^{\prime r})$ is well-defined.
For $L^{\prime }\in \mathcal{L}(\mathcal{C}^{\prime r}),$ $\widetilde{\theta 
}^{-1}(\mathcal{X}_{L^{\prime }})=\mathcal{X}_{\theta ^{-1}(L^{\prime })},$
i.e. $\mathbf{\theta }$ is continuous.

\item Suppose $\widetilde{\theta }(K_{1})=\widetilde{\theta }(K_{2})$ for
some $K_{1},K_{2}\in \mathrm{CPSpec}(\mathcal{C}^{r})$ with $K_{1}=\theta
^{-1}(K_{1}^{\prime }),$ $K_{2}=\theta ^{-1}(K_{2}^{\prime })$ where $%
K_{1}^{\prime },K_{2}^{\prime }\in \mathrm{CPSpec}(\mathcal{C}^{\prime r}).$
Then $K_{1}=\theta ^{-1}(K_{1}^{\prime })=\theta ^{-1}(\theta (\theta
^{-1}(K_{1}^{\prime })))=\theta ^{-1}(\theta (\theta ^{-1}(K_{1}^{\prime
})))=\theta ^{-1}(\theta (\theta ^{-1}(K_{2}^{\prime })))=\theta
^{-1}(K_{2}^{\prime })=K_{2}.$

\item By \cite[Proposition 4.7.]{Abu2006} $\mathrm{CPSpec}(\mathcal{C}^{r})=%
\mathcal{B}(\mathcal{C}^{r})\cap \mathrm{C}\mathrm{PSpec}(\mathcal{C}%
^{\prime r}),$ hence for $L\in \mathcal{L}(\mathcal{C}^{r})$ and $L^{\prime
}\in \mathcal{L}(\mathcal{C}^{\prime r})$ we have $\mathbf{\theta }^{-1}(%
\mathcal{V}_{L^{\prime }})=\mathcal{V}_{\theta ^{-1}(L^{\prime })},$ $%
\mathbf{\theta }(\mathcal{V}_{L})=\mathcal{V}_{\theta (L)}$ and $\mathbf{%
\theta }(\mathcal{X}_{L})=\mathcal{X}_{\theta (L)}.$

\item Since $\theta $ is an isomorphism, $\widetilde{\theta }$ is bijective
by \cite[Proposition 4.5.]{Abu2006}. In this case $\mathbf{\theta }$ and $%
\mathbf{\theta }^{-1}$ are obviously continuous (see \textquotedblleft
4\textquotedblright ).$\blacksquare $
\end{enumerate}
\end{Beweis}

\begin{ex}
\label{Ex1}(\cite[Example 1.1.]{NT2001}) Let $k$ be a field and $C:=k[X]$ be
the cocommutative $k$-coalgebra with $\Delta (X^{n}):=X^{n}\otimes _{k}X^{n}$
and $\varepsilon (X^{n}):=1$ for all $n\geq 0.$ For each $n\geq 0,$ set $%
C_{n}:=kX^{n}.$ Then $\mathrm{CPSpec}(C)=\mathcal{S}(C)=\{C_{n}\mid n\geq
0\}.$ Notice that

\begin{enumerate}
\item $\mathbf{Z}_{C}$ is discrete by Theorem \ref{T1}, hence $\mathbf{Z}%
_{C} $ is Lindelof (but not compact) by Theorem \ref{compact}.

\item $\mathrm{CPSpec}(C)$ is not connected: $\mathrm{CPSpec}(C)=\{C_{n}\mid
n\geq 1\}\cup \{C_{0}\}=\mathcal{X}_{\{k\}}\cup \mathcal{X}_{<X,X^{2},...>}$
(notice that $\mathrm{CPSpec}(C)$ is not subdirectly irreducible, compare
with Lemma \ref{it-irr}.
\end{enumerate}
\end{ex}

\begin{ex}
\label{Ex2}(\cite[Example 1.2.]{NT2001}) Let $k$ be a field and $C:=k[X]$ be
the cocommutative $k$-coalgebra with $\Delta
(X^{n}):=\dsum\nolimits_{j=1}^{n}X^{j}\otimes _{k}X^{n-j}$ and $\varepsilon
(X^{n}):=\delta _{n,0}$ for all $n\geq 0.$ For each $n\geq 0$ set $%
C_{n}:=<1,...,X^{n}>.$ For each $n\geq 1,$ $C_{n}\subseteq
(C_{n-1}:_{C}<kX^{n}>),$ hence not fully $C$-coprime and it follows that $%
\mathrm{CPSpec}(C)=\{k,C\}$ (since $k$ is simple, whence fully $C$-coprime
and $C^{\ast }\simeq k[[X]]$ is an integral domain, whence $C$ is fully
coprime). Notice that

\begin{enumerate}
\item $C$ is subdirectly irreducible with unique simple subcoalgebra $%
C_{0}=k;$

\item the converse of Remark \ref{simple-char} \textquotedblleft
3(d)\textquotedblright\ does not hold in general: $\mathrm{Corad}%
(C)=k\subseteq C_{1}$ while $\mathcal{X}_{C_{1}}=\{C\}\neq \varnothing $
(compare Proposition \ref{bireg} \textquotedblleft 2\textquotedblright ).

\item $\mathrm{CPSpec}(C)$ is connected, although $\mathcal{S}(C)\subsetneqq 
\mathrm{C}\mathrm{PSpec}(C)$ (see Lemma \ref{it-irr} \textquotedblleft
2\textquotedblright ).

\item $\mathbf{Z}_{C}$ is not $T_{1}$ by Theorem \ref{T1}, since $C\in 
\mathrm{CPSpec}(C)\backslash \mathcal{S}(C):$ in fact, if $C\in \mathcal{X}%
_{L_{1}}$ and $C_{0}\in \mathcal{X}_{L_{2}}$ for some $C$-subcoalgebras $%
L_{1},L_{2}\subseteq C,$ then $L_{2}=0$ (since $C$ is subdirectly
irreducible with unique simple subcoalgebra $C_{0}$); hence $\mathcal{X}%
_{L_{2}}=\{C_{0},C\}=\mathrm{CPSpec}(C)$ and $\mathcal{X}_{L_{1}}\cap 
\mathcal{X}_{L_{2}}=\mathcal{X}_{L_{1}}\neq \varnothing .$
\end{enumerate}
\end{ex}

\begin{remark}
As this paper extends results of \cite{NT2001}, several proofs and ideas are
along the lines of the original ones. However, our results are much more
general (as \cite{NT2001} is restricted to coalgebras over fields).
Moreover, we should warn the reader that in addition to the fact that
several results in that paper are redundant or repeated, several other
results are even \emph{absurd}, e.g. Proposition 2.8., Corollary 2.4. and
Theorem 2.4. (as noticed by Chen Hui-Xiang in his review; Zbl 1012.16041) in
addition to \cite[Lemma 2.6.]{NT2001} as we clarified in Remark \ref{PID}.
We corrected the statement of some of these results (e.g. Proposition \ref%
{bireg} corrects \cite[Lemma 2.6.]{NT2001}; while Proposition \ref{th-tel}
suggests a correction of \cite[Theorem 2.4.]{NT2001} which does not hold in
general as the counterexample \cite[5.20.]{Abu2006} shows). Moreover, we
improved some other results (e.g. applying Theorem \ref{T1} to coalgebras
over base fields improves and puts together several scattered results of 
\cite{NT2001}).
\end{remark}

\newpage

\textbf{Acknowledgments:} The author is grateful for the financial support
and the excellent research facilities provided by KFUPM.

\end{document}